\documentclass[12pt]{article}%
\usepackage{graphicx}
\usepackage[intlimits]{amsmath}
\usepackage{latexsym}
\usepackage{amsfonts}
\usepackage{amssymb}%
\makeatletter
\newfont{\footsc}{cmcsc10 at 8truept}
\newfont{\footbf}{cmbx10 at 8truept}
\newfont{\footrm}{cmr10 at 10truept}
\newtheorem{theorem}{Theorem}

\newtheorem{conjecture}[theorem]{Conjecture}
\newtheorem{corollary}[theorem]{Corollary}

\newtheorem{lemma}[theorem]{Lemma}

\newtheorem{proposition}[theorem]{Proposition}

\newenvironment{proof}[1][Proof]{\noindent{\textbf {#1}  }}  {\hfill$\Box$\bigskip}

\begin{document}

\title{Eigenvalues and degree deviation in graphs}
\author{Vladimir Nikiforov\\Department of Mathematical Sciences, University of Memphis, \\Memphis TN 38152, USA, email: \textit{vnkifrv@memphis.edu}}
\maketitle

\begin{abstract}
Let $G$ be a graph with $n$ vertices and $m$ edges and let $\mu\left(
G\right)  =\mu_{1}\left(  G\right)  \geq...\geq\mu_{n}\left(  G\right)  $ be
the eigenvalues of its adjacency matrix. Set $s\left(  G\right)  =\sum_{u\in
V\left(  G\right)  }\left\vert d\left(  u\right)  -2m/n\right\vert .$ We prove
that
\[
\frac{s^{2}\left(  G\right)  }{2n^{2}\sqrt{2m}}\leq\mu\left(  G\right)
-\frac{2m}{n}\leq\sqrt{s\left(  G\right)  }.
\]
In addition we derive similar inequalities for bipartite $G$.

We also prove that the inequality
\[
\mu_{k}\left(  G\right)  +\mu_{n-k+1}\left(  \overline{G}\right)
\geq-1-2\sqrt{2s\left(  G\right)  }%
\]
holds for every $k=1,...,n-1.$

We show that these inequalities are tight up to a constant factor.

Finally we prove that for every graph $G$ of order $n$,%
\[
\mu_{n}\left(  G\right)  +\mu_{n}\left(  \overline{G}\right)  \leq
-1-\frac{s^{2}\left(  G\right)  }{n^{3}}.
\]



\end{abstract}

\section{Introduction}

Our notation is standard (e.g., see \cite{Bol98}, \cite{CDS80}, and
\cite{HoJo88}); in particular, all graphs are defined on the vertex set
$\left\{  1,2,...,n\right\}  =\left[  n\right]  $ and $G\left(  n,m\right)  $
stands for a graph with $n$ vertices and $m$ edges. We write $\Gamma\left(
u\right)  $ for the set of neighbors of the vertex $u$ and set $d\left(
u\right)  =\left\vert \Gamma\left(  u\right)  \right\vert .$ Given a graph $G$
of order $n,$ we assume that the eigenvalues of the adjacency matrix of $G$
are ordered as $\mu\left(  G\right)  =\mu_{1}\left(  G\right)  \geq...\geq
\mu_{n}\left(  G\right)  $. As usual, $\overline{G}$ denotes the complement of
a graph $G.$

Collatz and Sinogowitz \cite{CoSi57} showed that $\mu\left(  G\right)
\geq2m/n$ for every graph $G=G\left(  n,m\right)  .$ Since equality holds if
and only if $G$ is regular, they proposed the value $\epsilon\left(  G\right)
=\mu\left(  G\right)  -2m/n$ as a relevant measure of irregularity of $G$. Two
other closely related measures of graph irregularity are the functions
\begin{align*}
var\left(  G\right)   &  =\frac{1}{n}\sum_{u\in V\left(  G\right)  }\left(
d\left(  u\right)  -\frac{2m}{n}\right)  ^{2},\\
s\left(  G\right)   &  =\sum_{u\in V\left(  G\right)  }\left\vert d\left(
u\right)  -\frac{2m}{n}\right\vert .
\end{align*}

Bell \cite{Bel92} compared $\epsilon\left(  G\right)  $ to $var\left(
G\right)  $ and showed that none of them could be preferred to the other one
as a measure of irregularity. He did not, however, give explicit inequalities
between $\epsilon\left(  G\right)  $ and $var\left(  G\right)  $. In this note
we prove that for every graph $G$ with $n$ vertices and $m$ edges,
\begin{equation}
\frac{var\left(  G\right)  }{2\sqrt{2m}}\leq\mu\left(  G\right)  -\frac{2m}%
{n}\leq\sqrt{s\left(  G\right)  } \label{mainin1}%
\end{equation}
Thus, in view of
\[
\frac{s^{2}\left(  G\right)  }{n^{2}}\leq var\left(  G\right)  \leq s\left(
G\right)  ,
\]
we also have
\[
\frac{s^{2}\left(  G\right)  }{2n^{2}\sqrt{2m}}\leq\mu\left(  G\right)
-\frac{2m}{n}\leq\sqrt{var\left(  G\right)  }.
\]
In addition we derive similar inequalities specifically for bipartite graphs.

Another well-known inequality involving graph eigenvalues is%
\begin{equation}
\mu_{k}\left(  G\right)  +\mu_{n-k+1}\left(  \overline{G}\right)  \leq-1,
\label{min2}%
\end{equation}
holding for every graph $G$ of order $n$ and every $k=1,...,n-1.$ Note that if
$G$ is regular, equality holds in (\ref{min2}) but the converse is not always
true (e.g., $G=K_{a,b},$ $b>a>2,$ $k=2)$. A natural problem is to find a lower
bound on $\mu_{k}\left(  G\right)  +\mu_{n-k+1}\left(  \overline{G}\right)  $
implying explicit equality in (\ref{min2}) for regular $G$. In this note we
show that for every $k=1,...,n-1,$
\begin{equation}
\mu_{k}\left(  G\right)  +\mu_{n-k+1}\left(  \overline{G}\right)
\geq-1-2\sqrt{2s\left(  G\right)  }. \label{mainin2}%
\end{equation}

We show that inequalities (\ref{mainin1}) and (\ref{mainin2}) are tight up to
a constant factor.

Finally we prove that for every graph $G$ of order $n$,%
\begin{equation}
\mu_{n}\left(  G\right)  +\mu_{n}\left(  \overline{G}\right)  \leq
-1-\frac{s^{2}\left(  G\right)  }{n^{3}},\label{mainin3}%
\end{equation}
implying that for any highly irregular graph $G$ either $\mu_{n}\left(
G\right)  $ or $\mu_{n}\left(  \overline{G}\right)  $ must be large in
absolute value.

Let us note that these results are readily applicable to the study of
quasirandom graph properties.

The rest of the note is organized as follows. In Section \ref{ER} we describe
algorithms for regularizing graphs with few edge changes. Section \ref{EVblo}
contains basic results about spectra of blown-up graphs. In Sections \ref{B1},
\ref{B2}, and \ref{B3} we prove inequalities (\ref{mainin1}), (\ref{mainin2}),
and (\ref{mainin3}).

\section{\label{ER}Efficient regularization}

Consider the following natural problem: given a graph $G,$ what is the minimum
number of edges $\rho\left(  G\right)  $ that must be changed to obtain a
regular graph. Writing $A\left(  G\right)  $ for the adjacency matrix of a
graph $G,$ we see that%
\[
\rho\left(  G\right)  =\frac{1}{2}\min\left\{  \left\Vert A\left(  G\right)
-A\left(  R\right)  \right\Vert _{2}:R\text{ is regular graph of order
}v\left(  G\right)  \right\}  .
\]
It is almost certain that the problem of estimating $\rho\left(  G\right)  $
has been raised and solved in the literature but, lacking a proper reference,
we shall solve it from scratch.

We first show that there exists a graph $R^{\ast}$ whose degrees differ by at
most one and such that
\[
\left\Vert A\left(  G\right)  -A\left(  R^{\ast}\right)  \right\Vert _{2}\leq
s\left(  G\right)  .
\]
Next we find a regular graph $R$ such that
\[
\left\Vert A\left(  R\right)  -A\left(  R^{\ast}\right)  \right\Vert _{2}<3n.
\]
Finally we show that for every graph $G,$%
\[
\rho\left(  G\right)  \geq s\left(  G\right)  /2
\]
implying that our upper bounds on $\rho\left(  G\right)  $ are not too far
from the best possible ones.

\subsection{Rough regularization}

The main result in this section is the following theorem.

\begin{theorem}
\label{th1}For every graph $G=G\left(  n,m\right)  ,$ there exists a graph
$R=G\left(  n,m\right)  $ such that $\Delta\left(  R\right)  \leq\delta\left(
R\right)  +1$ and $R$ differs from $G$ in at most $s\left(  G\right)  $ edges.
In particular, if $2m/n$ is integer then $R$ is $\left(  2m/n\right)  $-regular.
\end{theorem}

\begin{proof}
We shall describe a simple algorithm that produces the graph $R$ by deleting
and adding edges of $G.$ Set $d=\left\lfloor 2m/n\right\rfloor $.

\textbf{Step 1}

\emph{While }$\delta\left(  G\right)  <d$\emph{ and }$\Delta\left(  G\right)
>d+1$\emph{ select }$u,v$\emph{ with }$d\left(  u\right)  =\delta\left(
G\right)  $\emph{ and }$d\left(  v\right)  =\Delta\left(  G\right)  $\emph{.
Since }$\Gamma\left(  v\right)  \backslash\Gamma\left(  u\right)
\neq\varnothing,$\emph{ there exists }$w\in\Gamma\left(  v\right)
\backslash\Gamma\left(  u\right)  ;$\emph{ delete the edge }$vw$\emph{ and add
the edge }$uw.$

Write $G^{\prime}$ for the graph obtained upon exiting Step 1. Since Step 1 is
iterated as long as $\delta\left(  G\right)  <d$ and $\Delta\left(  G\right)
>d+1$, we have either $\delta\left(  G^{\prime}\right)  =d$ or $\Delta\left(
G^{\prime}\right)  =d+1;$ we may assume $\delta\left(  G^{\prime}\right)  =d,$
since the other case is reduced to this one by considering $\overline
{G^{\prime}}$.

If $\Delta\left(  G^{\prime}\right)  \leq d+1$ then terminate the procedure
with $R=G^{\prime}.$ Otherwise write $A$ for the set of vertices of degree
$d,$ $B$ for the set of vertices of degree $d+1,$ and $C$ for the set of
vertices of degree $d+2$ or higher.

\textbf{Step 2}

\emph{While }$C\neq\varnothing$\emph{ select }$u\in A,$\emph{ }$v\in C.$\emph{
Since }$\left\vert \Gamma\left(  v\right)  \right\vert >\left\vert
\Gamma\left(  u\right)  \right\vert ,$\emph{ we may select }$w\in\Gamma\left(
v\right)  \backslash\Gamma\left(  u\right)  ;$\emph{ delete the edge }%
$vw$\emph{ and add the edge }$uw.$

Write $R$ for the graph obtained after executing Step 2. Let $G^{\prime
},A,B,C$ be as defined prior to Step 2; set $\left\vert A\right\vert =k,$
$\left\vert C\right\vert =s$. Each iteration in Step 1 changes two edges and
decreases $s\left(  G\right)  $ by 2; therefore, after the execution of Step
1, at most $s\left(  G\right)  -s\left(  G^{\prime}\right)  $ edges of $G$ are
changed. Set
\[
l=\sum_{u\in C}\left(  d\left(  u\right)  -d-1\right)  .
\]
Each iteration in Step 2 changes two edges and decreases $l$ by 1; therefore,
there are $l$ iterations in Step 2 and at most $2l$ edges are changed. To
complete the proof we have to show that $S\left(  G^{\prime}\right)  \geq2l.$
From
\begin{align*}
2m/n  &  =\frac{1}{n}\sum_{u\in V\left(  G\right)  }d\left(  u\right)
=\frac{1}{n}\sum_{u\in A}d\left(  u\right)  +\frac{1}{n}\sum_{u\in B}d\left(
u\right)  +\frac{1}{n}\sum_{u\in C}d\left(  u\right) \\
&  =\frac{kd+\left(  n-k-s\right)  \left(  d+1\right)  +s\left(  d+1\right)
+l}{n}=d+\frac{n-k+l}{n}%
\end{align*}
it follows that $k>l.$ Furthermore,
\begin{align*}
s\left(  G^{\prime}\right)  -2l  &  =\sum_{u\in V}\left\vert d_{G^{\prime}%
}\left(  u\right)  -2m/n\right\vert -2l\\
&  =\sum_{u\in A}\left\vert d_{G^{\prime}}\left(  u\right)  -2m/n\right\vert
+\sum_{u\in B}\left\vert d_{G^{\prime}}\left(  u\right)  -2m/n\right\vert
+\sum_{u\in C}\left\vert d_{G^{\prime}}\left(  u\right)  -2m/n\right\vert
-2l\\
&  =k\frac{n-k+l}{n}+\left(  n-k-s\right)  \frac{k-l}{n}+s\frac{k-l}%
{n}-l=2\frac{\left(  k-l\right)  \left(  n-k\right)  }{n}>0,
\end{align*}
completing the proof.
\end{proof}

\subsubsection{Rough regularization of bipartite graphs}

Call a bipartite graph \emph{semiregular} if vertices belonging to the same
vertex class have equal degrees.

Let $G$ be a bipartite graph and $A,B$ be its vertex classes, $\left\vert
A\right\vert =a,$ $\left\vert B\right\vert =b.$ Define the function%
\[
s_{2}\left(  G\right)  =\sum_{u\in A}\left\vert d\left(  u\right)  -\frac
{m}{a}\right\vert +\sum_{u\in B}\left\vert d\left(  u\right)  -\frac{m}%
{b}\right\vert ;
\]
$s_{2}\left(  G\right)  $ is the equivalent to $s\left(  G\right)  $ for
bipartite graphs. Clearly, $s_{2}\left(  G\right)  =0$ if and only if $G$ is semiregular.

Modifying slightly the proof of Theorem \ref{th1} we obtain the following
special case for bipartite graphs.

\begin{theorem}
\label{th1.1}For every bipartite graph $G=G\left(  n,m\right)  $ with vertex
classes $A,B,$ there exists a bipartite graph $R=G\left(  n,m\right)  $ with
the same vertex classes such that:

\emph{(i) }$\left\vert d_{R}\left(  u\right)  -d_{R}\left(  v\right)
\right\vert \leq1$ for every $u,v$ belonging to the same vertex class;

\emph{(ii)} $R$ differs from $G$ in at most $s_{2}\left(  G\right)  $ edges.

In particular, if $m/\left\vert A\right\vert $ and $m/\left\vert B\right\vert
$ are integer then $R$ is semiregular.
\end{theorem}

\subsection{Fine regularization}

If we allow $m$ to change, we may further regularize the graph $R$ obtained in
Theorem \ref{th1}.

\begin{theorem}
\label{th2}Let the degrees of a graph $G=G\left(  n,m\right)  $ be either $d$
or $d+1.$ There exists an $r$-regular graph $R$ such that either $r=d$ or
$r=d+1,$ and $R$ differs from $G$ in at most $3n/2$ edges.
\end{theorem}

\begin{proof}
Write $A$ for the set of vertices of degree $d+1$ and $B$ for $V\left(
G\right)  \backslash A.$ Clearly either $\left\vert A\right\vert $ or
$\left\vert B\right\vert $ is even. We shall assume that $\left\vert
A\right\vert $ is even, otherwise we may apply the argument to the
complementary graph. Set $a=\left\vert A\right\vert .$ Our goal is to
construct a $d$-regular graph by changing at most $3a/2$ edges. We shall
describe a procedure constructing $R.$

\textbf{Step 1}

\emph{While }$E\left(  A\right)  \neq\varnothing,$\emph{ select }$uv\in
E\left(  A\right)  $\emph{ and remove it.}

\textbf{Step 2}.

\emph{While }$A\neq\varnothing,$\emph{ select two distinct }$u,v\in A$\emph{
and two disjoint vertices }$t\in\Gamma\left(  v\right)  ,$\emph{ }$w\in
\Gamma\left(  u\right)  .$\emph{ Delete the edges }$uw$ \emph{and} $vt;$\emph{
add the edge }$wt.$

The iteration in Step 2 may always be executed since, for every two distinct
$u,v\in A,$ there exist disjoint vertices $t\in\Gamma\left(  v\right)  $ and
$w\in\Gamma\left(  u\right)  .$ Indeed, if $\Gamma\left(  u\right)  \neq
\Gamma\left(  v\right)  ,$ select $w\in\Gamma\left(  u\right)  \backslash
\Gamma\left(  v\right)  .$ Since $d\left(  w\right)  =d<\left\vert
\Gamma\left(  v\right)  \right\vert ,$ there exists $t\in\Gamma\left(
v\right)  $ that is disjoint from $w$ and the assertion is proved. If
$\Gamma\left(  u\right)  =\Gamma\left(  v\right)  $ then $\Gamma\left(
u\right)  $ cannot induce a complete graph, since\ $\Gamma\left(  u\right)
\subset B$ and so the vertices in $\Gamma\left(  u\right)  $ have degree $d.$

Each iteration in Step 1 removes two vertices from $A$ and changes two edges.
Each iteration in Step 2 removes two vertices from $A$ and changes three
edges. Therefore, after changing at most $3\left\vert A\right\vert /2$ edges,
we obtain a $d$-regular graph $R$, as claimed.
\end{proof}

\subsection{Optimal regularization}

Summarizing Theorems \ref{th1} and \ref{th2}, we obtain the following corollary.

\begin{corollary}
\label{cor1}For every graph $G$ of order $n,$
\[
\rho\left(  G\right)  \leq s\left(  G\right)  +3n/2.
\]

\end{corollary}

It turns out that this bound is quite close to the optimal one, no matter what
the graph $G$ is. We shall show that
\[
\rho\left(  G\right)  \geq s\left(  G\right)  /2.
\]
Let $R$ be $r$-regular graph with $V\left(  R\right)  =V\left(  G\right)  .$
For every vertex $v\in V\left(  G\right)  ,$ we have%
\[
\left\vert \left(  \Gamma_{G}\left(  u\right)  \backslash\Gamma_{R}\left(
u\right)  \right)  \cup\left(  \Gamma_{R}\left(  u\right)  \backslash
\Gamma_{G}\left(  u\right)  \right)  \right\vert \geq d\left(  u\right)
+r-2\min\left(  d,r\right)  \geq\left\vert d\left(  u\right)  -r\right\vert .
\]
Hence, summing over all vertices $v\in V\left(  G\right)  $ we find that%
\[
2\rho\left(  G\right)  \geq\left\Vert A\left(  G\right)  -A\left(  R\right)
\right\Vert _{2}\geq\sum\left\vert d\left(  u\right)  -r\right\vert \geq
s\left(  G\right)  ,
\]
as claimed.

We note without proof that $\rho\left(  K_{a,b}\right)  \geq3s\left(
K_{a,b}\right)  /4$.

\section{\label{EVblo}The spectra of blown-up graphs}

In this section we introduce two operations on graphs and consider how they
affect graph spectra.

Let $G=G\left(  n,m\right)  $ and $t>0$ be integer. Write $G^{\left(
t\right)  }$ for the graph obtained by replacing each vertex $u\in V\left(
G\right)  $ by a set $V_{u}$ of $t$ vertices and joining $x\in V_{u}$ to $y\in
V_{v}$ if and only if $uv\in E\left(  G\right)  .$ Notice that $v\left(
G^{\left(  t\right)  }\right)  =tn.$ The following theorem holds.

\begin{theorem}
\label{thblo1}The\ eigenvalues of $G^{\left(  t\right)  }$ are $t\mu
_{1}\left(  G\right)  ,...,t\mu_{n}\left(  G\right)  $ together with $n\left(
t-1\right)  $ additional $0$'s.
\end{theorem}

Set $G^{\left[  t\right]  }=\overline{\overline{G}^{\left(  t\right)  }},$
i.e., $G^{\left[  t\right]  }$ is obtained from $G^{\left(  t\right)  }$ by
joining all vertices within $V_{u}$ for every $u\in V\left(  G\right)  ;$ note
also that $\overline{G^{\left(  t\right)  }}=\overline{G}^{\left[  t\right]
}.$ The following theorem holds.

\begin{theorem}
\label{thblo2}The\ eigenvalues of $G^{\left[  t\right]  }$ are $t\mu
_{1}\left(  G\right)  +t-1,...,t\mu_{n}\left(  G\right)  +t-1$ together with
$n\left(  t-1\right)  $ additional $\left(  -1\right)  $'s.
\end{theorem}

\section{\label{B1}Bounds on $\mu\left(  G\right)  $}

In this section we shall prove inequalities (\ref{mainin1}). Recall first the
inequality
\begin{equation}
\mu^{2}\left(  G\right)  \geq\frac{1}{n}\sum_{u\in V\left(  G\right)  }%
d^{2}\left(  u\right)  , \label{Hof1}%
\end{equation}
due to Hofmeister \cite{Hof88} and observe that Stanley's inequality
\cite{Sta87}%
\[
\mu\left(  G\right)  \leq-1/2+\sqrt{2m+1/4}%
\]
implies
\begin{equation}
\mu^{2}\left(  G\right)  \leq2m. \label{Sta1}%
\end{equation}

We thus find that
\begin{align*}
2\sqrt{2m}\left(  \mu\left(  G\right)  -2m/n\right)   &  \geq2\mu\left(
G\right)  \left(  \mu\left(  G\right)  -2m/n\right)  \geq\mu^{2}\left(
G\right)  -\left(  2m/n\right)  ^{2}\\
&  \geq\frac{1}{n}\sum_{u\in V\left(  G\right)  }d^{2}\left(  u\right)
-\left(  2m/n\right)  ^{2}=var\left(  G\right)  ,
\end{align*}
obtaining the lower bound in (\ref{mainin1}). To prove the upper bound we need
the following proposition.

\begin{proposition}
\label{pr1}If $G_{1}$ and $G_{2}$ are graphs with $V\left(  G_{1}\right)
=V\left(  G_{2}\right)  $ then
\[
\mu\left(  G_{1}\right)  -\mu\left(  G_{2}\right)  \leq\sqrt{2\left\vert
E\left(  G_{1}\right)  \backslash E\left(  G_{2}\right)  \right\vert }.
\]

\end{proposition}

\begin{proof}
Setting $G^{\prime}=\left(  V\left(  G_{1}\right)  ,E\left(  G_{1}\right)
\cup E\left(  G_{2}\right)  \right)  ,$ $G^{\prime\prime}=\left(  V\left(
G_{1}\right)  ,E\left(  G_{1}\right)  \backslash E\left(  G_{2}\right)
\right)  ,$ from Weyl's inequalities (\cite{HoJo88}, p. 181), we have
\[
\mu\left(  G_{1}\right)  \leq\mu\left(  G^{\prime}\right)  \leq\mu\left(
G_{2}\right)  +\mu\left(  G^{\prime\prime}\right)  .
\]
By (\ref{Sta1}), we have,%
\[
\mu\left(  G^{\prime\prime}\right)  \leq\sqrt{2\left\vert E\left(
G_{1}\right)  \backslash E\left(  G_{2}\right)  \right\vert },
\]
completing the proof.
\end{proof}

We shall deduce the upper bound in (\ref{mainin1}) essentially from Theorem
\ref{th1}.

\begin{theorem}
\label{th3}For every graph $G=G\left(  n,m\right)  ,$%
\[
\mu\left(  G\right)  -2m/n\leq\sqrt{s\left(  G\right)  }.
\]

\end{theorem}

\begin{proof}
Theorem \ref{th1} implies that there exists a graph $R=G\left(  n,m\right)  $
such that $\Delta\left(  R\right)  \leq\delta\left(  r\right)  +1$ and $R$
differs from $G$ in at most $s\left(  G\right)  $ edges. Since $e\left(
R\right)  =e\left(  G\right)  $ it follows that $\left\vert E\left(  G\right)
\backslash E\left(  R\right)  \right\vert =\left\vert E\left(  R\right)
\backslash E\left(  G\right)  \right\vert $ and so $2\left\vert E\left(
G\right)  \backslash E\left(  R\right)  \right\vert \leq s\left(  G\right)  .$
Hence, by Proposition \ref{pr1},
\begin{equation}
\mu\left(  G\right)  -2m/n\leq\mu\left(  G\right)  -\left\lceil
2m/n\right\rceil +1\leq\mu\left(  G\right)  -\mu\left(  R\right)
+1\leq1+\sqrt{s\left(  G\right)  }. \label{tempin1}%
\end{equation}

Notice that $v\left(  G^{\left(  t\right)  }\right)  =tn,$ $e\left(
G^{\left(  t\right)  }\right)  =t^{2}m,$ and $s\left(  G^{\left(  t\right)
}\right)  =t^{2}s\left(  G\right)  .$ Applying Theorem \ref{thblo1}, we also
see that
\[
\mu\left(  G^{\left(  t\right)  }\right)  =t\mu\left(  G\right)  .
\]
From (\ref{tempin1}) it follows that
\[
\left(  \mu\left(  G\right)  -2m/n\right)  t=\mu\left(  G^{\left(  t\right)
}\right)  -2e\left(  G^{\left(  t\right)  }\right)  /v\left(  G^{\left(
t\right)  }\right)  \leq1+\sqrt{s\left(  G^{\left(  t\right)  }\right)
}=1+t\sqrt{s\left(  G\right)  }.
\]
Hence, dividing by $t$ and letting $t$ tend to infinity, the desired
inequality follows.
\end{proof}

\subsection{Tightness of inequalities (\ref{mainin1})}

It is natural to ask how large $c$ could be so that the inequality%
\[
\mu\left(  G\right)  -\frac{2m}{n}\geq c\frac{s^{2}\left(  G\right)  }%
{n^{2}\sqrt{m}}%
\]
holds for every graph $G=G\left(  n,m\right)  .$ Taking the graph
$G=K_{n,n+1}$ for $n$ large enough, we see that $c$ may be at most $1/2$.

Similarly, let $c$ be such that the inequality%
\[
\mu\left(  G\right)  -2m/n\leq c\sqrt{s\left(  G\right)  }%
\]
holds for every graph $G=G\left(  n,m\right)  .$ Taking $G=K_{n}\cup K_{1}$ we
see that $c$ must be at least $1/\sqrt{2}.$

We venture the following conjecture.

\begin{conjecture}
For every graph $G$ of sufficiently large order $n$ and size $m$,,%
\[
\frac{s^{2}\left(  G\right)  }{2n^{2}\sqrt{m}}\leq\mu\left(  G\right)
-\frac{2m}{n}\leq\sqrt{s\left(  G\right)  /2}.
\]

\end{conjecture}

\subsection{Bounds on $\mu\left(  G\right)  $ when $G$ is bipartite}

It is possible to modify inequalities (\ref{mainin1}) to better suit bipartite graphs.

Let $G$ be a bipartite graph and $A,B$ be its vertex classes, $\left\vert
A\right\vert =a,$ $\left\vert B\right\vert =b.$ Then, by Rayleigh's principle
we have,%
\[
\mu\left(  G\right)  \geq e\left(  G\right)  /\sqrt{ab}.
\]
A careful analysis shows that equality is possible if and only if $G$ is
semiregular. In fact the following theorem holds.

\begin{theorem}
\label{th4}For every bipartite graph $G$ with vertex classes $A,B,$%
\[
\frac{s_{2}^{2}\left(  G\right)  }{2n^{2}\sqrt{\left\vert A\right\vert
\left\vert B\right\vert }}\leq\mu\left(  G\right)  -\frac{e\left(  G\right)
}{\sqrt{\left\vert A\right\vert \left\vert B\right\vert }}\leq\sqrt
{\frac{s_{2}\left(  G\right)  }{2}}.
\]

\end{theorem}

\begin{proof}
Let $\left\vert A\right\vert =a,$ $\left\vert B\right\vert =b,$ $e\left(
G\right)  =m,$ $v\left(  G\right)  =n.$ We start with the proof of the first
inequality. By the AM-QM inequality we have%
\begin{align*}
\sum_{u\in A}\left\vert d\left(  u\right)  -\frac{m}{a}\right\vert  &
\leq\sqrt{a\sum_{u\in A}\left(  d\left(  u\right)  -\frac{m}{a}\right)  ^{2}%
},\\
\sum_{u\in B}\left\vert d\left(  u\right)  -\frac{m}{b}\right\vert  &
\leq\sqrt{b\sum_{u\in B}\left(  d\left(  u\right)  -\frac{m}{b}\right)  ^{2}}.
\end{align*}
Hence, by Cauchy-Schwarz and inequality (\ref{Hof1}), we find that,%
\begin{align*}
s_{2}\left(  G\right)   &  \leq\sqrt{n}\sqrt{\sum_{u\in A}\left(  d\left(
u\right)  -\frac{m}{a}\right)  ^{2}+\sum_{u\in B}\left(  d\left(  u\right)
-\frac{m}{b}\right)  ^{2}}=\sqrt{n}\sqrt{\sum_{u\in V\left(  G\right)
}\left(  d^{2}\left(  u\right)  -\frac{m^{2}n}{ab}\right)  }\\
&  \leq n\sqrt{\mu^{2}\left(  G\right)  -\frac{m^{2}}{ab}}\leq n\sqrt{\left(
\mu\left(  G\right)  -\frac{m}{\sqrt{ab}}\right)  \left(  2\sqrt{ab}\right)
},
\end{align*}
proving the first inequality.

To prove the second inequality we first note the equivalent of Proposition
\ref{pr1} for bipartite graphs: if $G_{1}$ and $G_{2}$ are bipartite graphs
with the same vertex classes then
\[
\mu\left(  G_{1}\right)  -\mu\left(  G_{2}\right)  \leq\sqrt{\left\vert
E\left(  G_{1}\right)  \backslash E\left(  G_{2}\right)  \right\vert }.
\]
Note that the coefficient $2$ under the square root is missing here, since
$\mu\left(  G\right)  \leq\sqrt{e\left(  G\right)  }$ for bipartite $G$
(Cvetkovi\'{c} \cite{Cve72}, also \cite{CDS80}, p. 92 Theorem 3.19).

Theorem \ref{th1.1} implies that there exists a graph $R=G\left(  n,m\right)
$ with vertex classes $A,B$ such that $\left\vert d_{R}\left(  u\right)
-d_{R}\left(  v\right)  \right\vert \leq1$ for every $u,v$ belonging to the
same vertex class and $R$ differs from $G$ in at most $s_{2}\left(  G\right)
$ edges. Since $e\left(  R\right)  =e\left(  G\right)  $ it follows that
$\left\vert E\left(  G\right)  \backslash E\left(  R\right)  \right\vert
=\left\vert E\left(  R\right)  \backslash E\left(  G\right)  \right\vert $ and
so $2\left\vert E\left(  G\right)  \backslash E\left(  R\right)  \right\vert
\leq s_{2}\left(  G\right)  .$ Hence, by Proposition \ref{pr1},
\[
\mu\left(  G\right)  -\mu\left(  R\right)  \leq\sqrt{\frac{s_{2}\left(
G\right)  }{2}}.
\]
Applying the inequality $\mu\left(  G\right)  \leq\max_{uv\in E\left(
G\right)  }\sqrt{d\left(  u\right)  d\left(  v\right)  },$ due to Berman and
Zhang \cite{BeZh01}, we find that
\[
\mu\left(  R\right)  \leq\sqrt{\left(  \frac{m}{a}+1\right)  \left(  \frac
{m}{b}+1\right)  }\leq\sqrt{\frac{m^{2}}{ab}+\frac{mn}{ab}+1}<\frac{m}%
{\sqrt{ab}}+\sqrt{n+1}%
\]
and so,%
\[
\mu\left(  G\right)  -\frac{m}{\sqrt{ab}}\leq\sqrt{\frac{s_{2}\left(
G\right)  }{2}}+\sqrt{n+1}.
\]
Now, applying the final argument from the proof of Theorem \ref{th3}, the
desired inequality follows.
\end{proof}

\section{\label{B2}A lower bound on $\mu_{k}\left(  G\right)  +\mu
_{n-k+1}\left(  \overline{G}\right)  $}

The main goal of this section is the proof of inequality (\ref{min2}). By
Weyl's inequalities (\cite{HoJo88}, p. 181), for every graph $G$ of order $n,$
we have
\[
\mu_{k}\left(  G\right)  +\mu_{n-k+1}\left(  \overline{G}\right)  \leq\mu
_{k}\left(  K_{n}\right)  =-1.
\]

\begin{theorem}
For every $k=1,...,n-1$%
\[
\mu_{k}\left(  G\right)  +\mu_{n-k+1}\left(  \overline{G}\right)
\geq-1-2\sqrt{2s\left(  G\right)  }%
\]

\end{theorem}

\begin{proof}
By Corollary \ref{cor1} there exists a regular graph $R$ that differs from $G$
in at most $s\left(  G\right)  +3n/2$ edges. Then, by Weyl's inequalities,
\begin{align*}
\mu_{k}\left(  A\left(  G\right)  \right)  +\mu_{1}\left(  A\left(  R\right)
-A\left(  G\right)  \right)   &  \geq\mu_{k}\left(  A\left(  R\right)
\right)  ,\\
\mu_{n-k+1}\left(  A\left(  \overline{G}\right)  \right)  +\mu_{1}\left(
A\left(  \overline{R}\right)  -A\left(  \overline{G}\right)  \right)   &
\geq\mu_{n-k+1}\left(  A\left(  \overline{R}\right)  \right)  .
\end{align*}
Furthermore, by
\begin{align*}
\mu_{1}\left(  A\left(  R\right)  -A\left(  G\right)  \right)   &  \leq
\sqrt{\left\Vert A\left(  R\right)  -A\left(  G\right)  \right\Vert _{2}%
}=\sqrt{2s\left(  G\right)  +3n}\\
\mu_{1}\left(  A\left(  \overline{R}\right)  -A\left(  \overline{G}\right)
\right)   &  \leq\sqrt{\left\Vert A\left(  \overline{R}\right)  -A\left(
\overline{G}\right)  \right\Vert _{2}}=\sqrt{2s\left(  G\right)  +3n},
\end{align*}
we find that
\begin{align*}
\mu_{k}\left(  G\right)  +\mu_{n-k+1}\left(  \overline{G}\right)   &  \geq
\mu_{k}\left(  A\left(  R\right)  \right)  +\mu_{n-k+1}\left(  A\left(
\overline{R}\right)  \right)  -2\sqrt{2s\left(  G\right)  +3n}\\
&  =-1-2\sqrt{2s\left(  G\right)  +3n}.
\end{align*}

Suppose now that $t$ is sufficiently large and consider the graphs $G^{\left(
t\right)  }$ and $\overline{G^{\left(  t\right)  }}.$ By Theorem \ref{thblo1}
we have%
\[
\mu_{k}\left(  G^{\left(  t\right)  }\right)  =t\mu_{k}\left(  G\right)  .
\]
Similarly in view of and $\overline{G^{\left(  t\right)  }}=\overline
{G}^{\left[  t\right]  }$ and Theorem \ref{thblo2},%
\[
\mu_{nt-k+1}\left(  \overline{G^{\left(  t\right)  }}\right)  \leq\min\left\{
t\mu_{n-k+1}\left(  \overline{G}\right)  +t-1,-1\right\}
\]
Since, $s\left(  G^{\left(  t\right)  }\right)  =t^{2}s\left(  G\right)  ,$ we
see that%
\begin{align*}
t\mu_{k}\left(  G\right)  +t\mu_{n-k+1}\left(  \overline{G}\right)   &
\geq\mu_{k}\left(  G^{\left(  t\right)  }\right)  +\mu_{k}\left(
\overline{G^{\left(  t\right)  }}\right)  -t+1\geq-t-2\sqrt{2s\left(
G^{\left(  t\right)  }\right)  +3nt}\\
&  =-t-2t\sqrt{2s\left(  G\right)  +3n/t}.
\end{align*}
Dividing by $t$ and letting $t$ tend to infinity, we obtain the desired inequality.
\end{proof}

For the graph $G=K_{1,n}$ we have $s\left(  G\right)  =2\frac{n-1}{n+1}$ and
$\mu_{n+1}\left(  G\right)  +\mu_{2}\left(  \overline{G}\right)  =-1-\sqrt
{n}.$ Hence,%
\[
\mu_{n+1}\left(  G\right)  +\mu_{2}\left(  \overline{G}\right)  =-1-\left(
\frac{1}{\sqrt{2}}+o\left(  1\right)  \right)  \sqrt{s\left(  G\right)  },
\]
implying that inequality (\ref{mainin2}) is tight up to a constant factor less
than $4.$

\section{\label{B3}An upper bound on $\mu_{n}\left(  G\right)  +\mu_{n}\left(
\overline{G}\right)  $}

The main result in this section is the proof of inequality (\ref{mainin3}). We
start with an auxiliary result.

\begin{lemma}
\label{lear} For every graph $G$ of order $n$ there exists an $\left\lfloor
n/2\right\rfloor $-set $S\subset V\left(  G\right)  $ such that
\[
e\left(  V\left(  G\right)  \backslash S\right)  -e\left(  S\right)  \geq
\frac{1}{2}s\left(  G\right)  .
\]

\end{lemma}

\begin{proof}
Note first that for any $a$ we have
\[
\sum_{i=1}^{n}\left\vert d_{i}-a\right\vert \geq\sum_{i=1}^{n}\left\vert
d_{i}-\frac{2m}{n}\right\vert =s\left(  G\right)  .
\]
Let $d\left(  1\right)  \leq d\left(  2\right)  \leq....\leq d\left(
n\right)  $ be the degree sequence of $G$ and set $V=\left[  n\right]  .$ For
every $1\leq k\leq n,$ letting $S=[k],$ we have
\[
\sum_{u\in V\backslash S}d\left(  u\right)  -\sum_{u\in S}d\left(  u\right)
=2e\left(  S\right)  +e\left(  S,V\backslash S\right)  -2e\left(  V\backslash
S\right)  -e\left(  S,V\backslash S\right)  =2e\left(  V\backslash S\right)
-2e\left(  S\right)  .
\]

Assume first $n$ even, $n=2k$. Letting $a=\left(  d\left(  k\right)  +d\left(
k+1\right)  \right)  /2$ and $S=[k],$ we have
\[
\sum_{u\in V\backslash S}d\left(  u\right)  -\sum_{u\in S}d\left(  u\right)
=\sum_{u\in V\backslash S}\left(  d\left(  u\right)  -a\right)  +\sum_{u\in
S}\left(  a-d\left(  u\right)  \right)  =\sum_{u\in V}\left\vert
d_{i}-a\right\vert \geq s\left(  G\right)  ,
\]
proving the assertion for even $n.$ 

Let now $n$ be odd, $n=2k+1$. Letting $a=d_{k+1}$ and $S=[k],$ we have
\[
\sum_{u\in V\backslash S}d\left(  u\right)  -\sum_{u\in S}d\left(  u\right)
=\sum_{u\in V\backslash S}\left(  d\left(  u\right)  -a\right)  +\sum_{u\in
S}\left(  a-d\left(  u\right)  \right)  =\sum_{u\in V}\left\vert d\left(
u\right)  -a\right\vert \geq s\left(  G\right)  ,
\]
proving the assertion for odd $n$ as well.
\end{proof}

\begin{theorem}
For every graph $G$ of order $n,$%
\[
\mu_{n}\left(  G\right)  +\mu_{n}\left(  \overline{G}\right)  \leq
-1-\frac{s^{2}\left(  G\right)  }{n^{3}}.
\]

\end{theorem}

\begin{proof}
From the interlacing theorem of Haemers (see, e.g., \cite{Ha95},
\cite{BoNi04}), for every bipartition of $V\left(  G\right)  =V_{1}\cup V_{2}$
we have%
\begin{equation}
\mu_{n}\left(  G\right)  \leq\frac{e\left(  V_{1}\right)  }{\left\vert
V_{1}\right\vert }+\frac{e\left(  V_{2}\right)  }{\left\vert V_{2}\right\vert
}-\sqrt{\left(  \frac{e\left(  V_{1}\right)  }{\left\vert V_{1}\right\vert
}-\frac{e\left(  V_{2}\right)  }{\left\vert V_{2}\right\vert }\right)
^{2}+\frac{e\left(  V_{1},V_{2}\right)  ^{2}}{\left\vert V_{1}\right\vert
\left\vert V_{2}\right\vert }}.\label{Hamin}%
\end{equation}

Assume $n$ even and let $V\left(  G\right)  =V_{1}\cup V_{2}$ be a bipartition
such that $\left\vert V_{1}\right\vert =\left\vert V_{2}\right\vert =n/2,$ and
$e\left(  V_{1}\right)  -e\left(  V_{2}\right)  \geq s\left(  G\right)  /2.$
Letting $e_{1}=e\left(  V_{1}\right)  $, $e_{2}=e\left(  V_{2}\right)  ,$
$e_{3}=e\left(  V_{1},V_{2}\right)  ,$ $s=s\left(  G\right)  ,$ from
(\ref{Hamin}), after some simple algebra, we obtain
\begin{equation}
\frac{n}{2}\mu_{n}\left(  G\right)  \leq e_{1}+e_{2}-\sqrt{\left(  e_{1}%
-e_{2}\right)  ^{2}+e_{3}^{2}}\leq e_{1}+e_{2}-\sqrt{\frac{s^{2}}{4}+e_{3}%
^{2}}.\label{Hamin1}%
\end{equation}

Note that $s\left(  G\right)  <n^{2}$ and $e\left(  V_{1},V_{2}\right)  \leq
n^{2}/4;$ thus, we have
\[
\frac{s^{4}}{9n^{4}}+\frac{2e_{3}s^{2}}{3n^{2}}+e_{3}^{2}\leq s^{2}\left(
\frac{1}{9}+\frac{1}{6}\right)  +e_{3}^{2}\leq\frac{s^{2}}{4}+e_{3}^{2},\text{
}%
\]
and so,
\[
\sqrt{\frac{s^{2}}{4}+e_{3}^{2}}\geq\frac{s^{2}}{3n^{2}}+e_{3}.
\]
Hence, from (\ref{Hamin1}), it follows that
\[
\frac{n}{2}\mu_{n}\left(  G\right)  \leq e_{1}+e_{2}-e_{3}-\frac{s^{2}}%
{3n^{2}}.
\]

Since $s\left(  G\right)  =s\left(  \overline{G}\right)  ,$ we see also that%
\[
\frac{n}{2}\mu_{n}\left(  \overline{G}\right)  \leq\binom{n/2}{2}-e_{1}%
+\binom{n/2}{2}-e_{2}-\frac{n^{2}}{4}+e_{3}-\frac{s^{2}}{3n^{2}}%
\]
and hence,
\[
\frac{n}{2}\left(  \mu_{n}\left(  G\right)  +\mu_{n}\left(  \overline
{G}\right)  \right)  \leq2\binom{n/2}{2}-\frac{n^{2}}{4}-\frac{2s^{2}}{3n^{2}%
}=-\frac{n}{2}-\frac{2s^{2}}{3n^{2}},
\]
proving the assertion for even $n.$

To prove the assertion for odd $n$ observe that if $t$ is even, for the graph
$G^{\left(  t\right)  }$ we have%
\[
t\mu_{n}\left(  G\right)  +t\mu_{n}\left(  \overline{G}\right)  +t-1=\mu
_{tn}\left(  G^{\left(  t\right)  }\right)  +\mu_{tn}\left(  \overline
{G^{\left(  t\right)  }}\right)  \leq-1-\frac{t^{4}s^{2}}{t^{3}n^{3}}.
\]
Dividing by $t$ and letting $t$ tend to infinity, the assertion follows for
odd $n$ as well.
\end{proof}

\end{document}